\documentclass[a4paper,11pt]{article}
\usepackage[psamsfonts]{amssymb}
\usepackage[psamsfonts]{eucal}

\usepackage{amsmath}
\usepackage{theorem}
\usepackage{color}
\newtheorem{theorem}{Theorem}
\newtheorem{lemma}{Lemma}

\newtheorem{prop}{Proposition}

\begin{document}

\title{ Complex surfaces and null conformal Killing vector fields}

\author{J. Davidov, G. Grantcharov, O. Mushkarov}

\date{}

\maketitle

\begin{abstract}

We study the relation between the existence of null conformal
Killing vector fields and existence of compatible complex and
para-hypercomplex structures on a pseudo-Riemannian manifold with
metric of signature $(2,2)$. We establish first the topological
types of pseudo-Hermitian surfaces admitting a nowhere vanishing
null vector field. Then we show that a pair of orthogonal, pointwise
linearly independent, null, conformal Killing vector fields defines
a para-hyperhermitian structure and use this fact for a
classification of the smooth compact four-manifolds admitting such  a pair
of vector fields. We also provide examples of neutral  metrics with
two orthogonal, pointwise linearly independent, null Killing vector
fields on most of  these manifolds.

\end{abstract}

\section{Introduction}

 The notion of self-duality for metrics of split
signature $(2,2)$  called also neutral metrics, has been studied for
a long time and related to various fields like integrable systems
\cite{Ward}  and superstring with $N=2$ supersymmetry \cite{OV1,
OV2}. Many results  can directly be extended from the positive
definite to the split signature case but there are also some
significant differences. Most of the similarities , as well as the
differences, are based on the fact that the metrics of signature
$(2,2)$ are related to the split-quaternions
$$
\mathbb{H}' =\{q = a+bi +cs + dt \in \mathbb{R}^4 | s^2=t^2=1=-i^2,
t = is, is=-si \}
$$
in the same way as  the positive definite metrics to  the
quaternions . For example, the similarities include a twistor space
construction  \cite{BDM} and a spinor approach  \cite{DW1}, while
the differences come from the fact that the self-duality equations
in split signature are not elliptic but ultra-hyperbolic. Most of
the research has been focused on local study in view of the fact
that for ultra-hyperbolic equations the global properties are harder
to understand.

There has also been some interest in compact 4-manifolds admitting
neutral metrics and compatible complex or para-hypercomplex
structures since topological information like the Kodaira
classification of compact complex surfaces allows one to study
global properties. Important results in this direction are the
classifications of compact pseudo-K\"ahler Einstein and
para-hyperk\"ahler surfaces obtained by Petean \cite{P} and Kamada
\cite{Kamada02}, respectively. These structures  appear in
 {\cite{OV2}} as models for superstring theory with
N=2 supersymmetry and  \cite{Hit} in relation to deformation spaces
of harmonic maps from Riemann surfaces into Lie groups. Note also
that such structures have been used recently by B. Klingler
\cite{Klin} in his proof of the Chern conjecture for affine
manifolds.

In our previous papers \cite{DGMY09, DGMY} we initiated the study of
compact para-hyperhermitian surfaces, which are the neutral analog
of the hyperhermitian 4-manifolds. These surfaces are anti-self-dual
as in the positive definite case, but  in contrast to the well-known
classification of compact hyperhermitian surfaces \cite{Boyer},
there are many more compact examples of para-hyperhermitian
surfaces. In our study, we noticed also that any two pointwise
independent, null, orthogonal and parallel vector fields on a
4-manifold with a neutral metric define a para-hyperk\"aher
structure and conversely, any compact para-hyperk\"aher surface
admits two vector fields with the above properties \cite{Kamada02}.
In this note, we weaken these conditions and study the neutral
Hermitian surfaces admitting non-vanishing null conformal Killing
vector fields. We observe first that the existence of a
non-vanishing null vector field on a compact neutral Hermitian
surface leads to a topological restriction which implies a rough
classification of these surfaces. Next, we show that the existence
of two point-wise independent, null and orthogonal conformal Killing
vector fields leads to the existence of a para-hyperhermitian
structure. These surfaces were studied in \cite{DGMY09, DGMY} and
here we consider the problem for existence of two Killing vector
fields having the above properties with respect to the corresponding
metrics.

Now we describe shortly the content of the paper. In Section 2, some
known facts about  almost para-hyperhermitian structures are
collected and in Lemma 2  they are characterized in terms of
2-forms. Then in Section 3 we study the complex surfaces admitting
neutral Hermitian metrics with non-vanishing null  vector fields.
After proving some useful local properties of these surfaces, in
Theorem~\ref{compact-1-field}, we provide  a classification in the
compact case. In Theorem~\ref{RH-ASD} we notice a property analogous
to the positive definite case - if a neutral Hermitian surface
admits  complete null Killing vector field which is not real
holomorphic, then the metric is anti-self-dual. Section 4 is devoted
to the study of the neutral 4-manifolds admitting two orthogonal and
null conformal Killing vector fields which are linearly independent
at each point. We show that any such a manifold admits a
para-hypercomplex structure which is compatible with the given
metric (Theorem~\ref{thm1}) and  provide a classification of these
manifolds in the compact case (Theorem~\ref{existence of 2 fields}).
In Theorem~\ref{existence of 2 fields}, we also exhibit explicit
examples of neutral metrics having two Killing vector fields with
the required properties on most of the possible smooth
$4$-manifolds.

\vspace{.1in}

{\bf Acknowledgments}.The first and the third named authors are partially supported by the Bulgarian National Science Fund, Ministryof Education and Science of Bulgaria under contract KP-06-N52/3.

The second named author is supported by the Simons Foundation
Collaboration grant number 853269. He wants to express his
gratitude to the Institute of Mathematics and Informatics of the
Bulgarian Academy of Sciences, Sofia for the warm hospitality during
his visits to work on the project.

\section{Para-hyperhermitian structures}

Recall that an {\it almost para-hypercomplex} structure on a smooth
manifold $M$ is a triple $(I,S,T)$ of anti-commuting endomorphisms
of  its tangent bundle $TM$ with $I^2=-Id$ and $S^2=T^2=Id$, $T=IS$.
Such a structure is called {\it para-hypercomplex} if $I,S,T$
satisfy the integrability condition $N_I=N_S=N_T=0$, where
$$
N_A(X,Y)=A^2[X,Y]+[AX,AY]-A[AX,Y]-A[X ,AY]
$$
is the Nijenhuis tensor of $A=I,S,T$. Note that if two of the
structures $I,S,T$ are integrable, the third one is also integrable
\cite{Kamada02}.

If there exists a pseudo-Riemannian metric $g$ for which the
endomorphisms $I,S,T$ are skew-symmetric,  $(g,I,S,T)$ is called an
{\it almost para-hyperhermitian} structure. Such a metric is
necessarily neutral, i.e., has split signature. We also say that the
metric $g$ is compatible with the almost para-hypercomplex structure
$(I,S,T)$ as well as that $(I,S,T)$ is compatible with $g$. Every
almost para-hypercomplex structure on a 4-manifold locally admits a
compatible metric but a globally defined one may not exist. More
precisely, such a structure determines a conformal class of
compatible metrics up to a double cover of $M$ \cite{DGMY09}.
Examples of almost para-hypercomplex structures that do not admit
compatible metrics are given in \cite{DGMY09, DGMY}.

Let $(g,I,S,T)$ be an almost para-hyperhermitian structure on a
$4$-manifold. Then  one can define three fundamental $2$-forms
$\Omega_i$, $i=1,2,3$, setting $$\Omega_1(X,Y)=g(IX,Y),
\hspace{.1in} \Omega_2(X,Y)=g(SX,Y), \hspace{.1in}
\Omega_3(X,Y)=g(TX,Y).$$ Note that the form
$\Omega=\Omega_2+i\Omega_3$ is of type $(2,0)$ with respect to $I$.
As in the definite case, the corresponding Lee forms are defined by
$$\theta_1= -\delta\Omega_1\circ I, \hspace{.1in}\theta_2=-\delta\Omega_2\circ S,
\hspace{.1in}\theta_3=-\delta\Omega_3\circ T,$$ where $\delta$ is
the co-differential with respect to $g$.

It is well-known \cite{Boyer, GT, Iv-Zam, Kamada02} that the
structures $I,S,T$ are intergrable if and only if
$\theta_1=\theta_2=\theta_3$. Thus, for a para-hyperhermitian
structure, we have just one Lee form $\theta$; it satisfies the
identities
$$d\Omega_i=\theta\wedge\Omega_i,~ i=1,2,3.$$ When additionally
the three 2-forms $\Omega_i$ are closed, i.e. $\theta=0$, the
para-hyperhermitian structure is  {\it para-hyperk\"ahler} ( also
hypersympectic or neutral hyperk\"ahler). When $d\theta=0$ the
structure is {\it locally conformally para-hyperk\"ahler}. { We note
that, in dimension 4, the para-hyperhermitian metrics are self-dual
and the para-hyperk\"ahler metrics are self-dual and Ricci-flat
\cite{Kamada02}. It is well-known that every hyperhermitian
structure on a 4-dimensional compact manifold is locally conformally
hyperk\"ahler \cite{Boyer}, but it has been shown in \cite[Theorem
9]{DGMY}, that this is not true in the indefinite case.

An almost para-hyperhermitian structure on a $4$-manifold can be
characterized by means of the forms $\Omega_i$ in the following way
\cite{Kamada02}.
\begin{prop}\label{phe}
Every almost para-hyperhermitian structure on a $4$-manifold is
uniquely determined by three non-degenerate $2$-forms $(\Omega_1,
\Omega_2, \Omega_3)$  such that
$$
-\Omega_1^2 = \Omega_2^2 = \Omega_3^2, \quad \Omega_l \wedge\Omega_m
= 0, \, 1\leq l \neq m \leq 3.
$$
This structure is para-hyperhermitian  if and only if there is a
$1$-form $\theta$ such that $d\Omega_l=\theta\wedge\Omega_l$,
$l=1,2,3$,
\end{prop}
 Note that given $2$-forms $(\Omega_1,
\Omega_2, \Omega_3)$ with the above properties the endomorphisms
$I,S,T$ and the metric $g$ of the almost para-hyperhermitian
structure they determine are defined by
$$
\Omega_3(IX,Y)=\Omega_2(X,Y),\quad \Omega_1(SX,Y)=\Omega_3(X,Y),
\quad \Omega_2(TX,Y)=-\Omega_1(X,Y),
$$
and
$$
 g(X,Y)=\Omega_1(X,IY),\quad X,Y\in TM.
$$
One can show that
$$
\Omega_1(X,IY)=-\Omega_2(X,SY)=-\Omega_3(X,TY),
$$
so $\Omega_1$, $\Omega_2$, $\Omega_3$ are the fundamental $2$-forms
of $I$, $S$, $T$.

\smallskip

We use later on the following characterization of the almost
para-hyperhermitian structures with a given almost complex
structure.

\begin{lemma}\label{APHHS}
An almost complex $4$-manifold $(M,I)$ admits an almost
para-hyperhermitian structure with almost complex structure $I$ if
and only if there is a $(2,0)$-form $\Omega$ and a non-degenerate
real $2$-form $\omega$ such that $\Omega\wedge\overline\Omega=
-2\omega^2$ and $\Omega\wedge\omega=0$. If $\Omega$ and $\omega$
satisfy these conditions, $\omega$, $Re\,\Omega$, $Im\,\Omega$ are
the fundamental 2-forms of the almost para-hyperhermitian structure.

\end{lemma}

\noindent {\bf Proof}.  Suppose we are given $2$-forms $\Omega$ and
$\omega$ satisfying the conditions of the lemma. Set
$\Omega_1=\omega$, $\Omega_2=Re\,\Omega$, $\Omega_3=Im\,\Omega$. We
have $\Omega\wedge\Omega=0$ since $\Omega$ is a $(2,0)$-form
 and $dim\,M=4$. This implies $\Omega_2^2=\Omega_3^2$. Thus,
$2\Omega_2^2=2\Omega_3^2=\Omega\wedge\overline\Omega= -2\omega^2\neq
0$. Hence the $2$-forms $\Omega_1, \Omega_2,\Omega_3$ are
non-degenerate. The identities $\Omega\wedge\Omega=0$ and
$\Omega\wedge\overline\Omega=-2\omega^2$ give
$-\Omega_1^2=\Omega_2^2=\Omega_3^2$ and $\Omega_2\wedge \Omega_3=0$.
Moreover, the identity $\Omega\wedge\omega = 0$ is equivalent to
$\Omega_2\wedge\Omega_1=\Omega_3\wedge\Omega_1=0$. By
Proposition~\ref{phe}, there exists an almost para-hyperhermitian
structure $(g,I',S,T)$ such that $\Omega_1,\Omega_2,\Omega_3$ are
the respective fundamental 2-forms. The form
$\Omega=\Omega_2+i\Omega_3$ is of type $(2,0)$ with respect to both
$I$ and $I'$. It follows that $\Omega_2(X,Y)=\Omega_3(X,IY)$ and
$\Omega_2(X,Y)=\Omega_3(X,I'Y)$. Hence $I=I'$ since the form
$\Omega_3(X,Y)=g(X,TY)$ is non-degenerate.

The rest part of the lemma is obvious.

\hfill$\Box$

\noindent {\bf Remark}. If $\Omega$ and $\omega$ satisfy the
conditions of Lemma~\ref{APHHS}, then $\omega$ is of type $(1,1)$
w.r.t. $I$ since $\omega(X,Y)=g(IX,Y)$ for a
 neutral metric $g$ compatible with the almost complex structure $I$.

 \medskip


\

\section{Neutral Hermitian surfaces with a null vector field}

We begin with preliminary algebraic considerations.

Let $M$ be an oriented $4$-manifold with a neutral metric $g$ and
let $(E_1,E_2,E_3,E_4)$ be an oriented orthonormal frame such that
$g(E_1,E_1)=g(E_2,E_2)=1$ and $g(E_3,E_3)=g(E_4,E_4)=-1$. The Hodge
star operator of $g$ acts as an involution  on the bundle of
$2$-vectors $\Lambda^2 TM$ and is given by
$$\ast E_1\wedge E_2=E_3\wedge E_4,\quad
\ast E_1\wedge E_3=E_2\wedge E_4,\quad \ast E_1\wedge E_4=-E_2\wedge
E_3.$$ Let $\Lambda^2_{\pm}TM$ be the subbundles of $\Lambda^2TM$
corresponding the the eigenvalues $\pm 1$ of the Hodge operator. Set
\begin{equation}\label{basis}
\begin{array}{c}
s_1^{+}=E_1\wedge E_2+E_3\wedge E_4,\quad\quad
s_1^{-}=E_1\wedge E_2-E_3\wedge E_4,\\[6pt]
s_2^{+}=E_1\wedge E_3-E_4\wedge E_2,\quad\quad
s_2^{-}=E_1\wedge E_3+E_4\wedge E_2,\\[6pt]
s_3^{+}=E_1\wedge E_4-E_2\wedge E_3,\quad\quad s_3^{-}=E_1\wedge
E_4+E_2\wedge E_3.
\end{array}
\end{equation}
Then $\{s_1^{\pm},s_2^{\pm},s_3^{\pm}\}$ is an orthonormal frame of
$\Lambda^2_{\pm}TM$.

\smallskip

Recall that a $2$-plane in $TM$ is called self-dual or
$\alpha$-plane if there is a basis $a, b$ of the plane such that
$\ast a\wedge b=a\wedge b$. This condition does not depend on the
choice of the basis $a,b$. Similarly, a $2$-plane in $TM$ is called
anti-self-dual, or $\beta$-plane if there is a basis $a,b$ with
$\ast a\wedge b=-a\wedge b$.

If $\Pi$ is an isotropic $2$-plane in a tangent space $T_pM$, one
can find a complementary isotropic $2$-plane $\Pi'$ and bases
$a_1,a_2$ of $\Pi$ and  $a_1', a_2'$ of $\Pi'$ such that
$g(a_i,a_j')=\frac{1}{2}\delta_{ij}$. Then $e_1=a_1+a_1'$,
$e_2=a_2+a_2'$, $e_3=a_1-a_1'$, $e_4=a_2-a_2'$ is an orthonormal
basis of $T_pM$ such that $||e_1||^2=||e_2||^2=1$,
$||e_3||^2=||e_4||^2=-1$. We have $4 a_1\wedge a_2=(e_1+e_3)\wedge
(e_2+e_4)=(e_1\wedge e_2+e_3\wedge e_4)+(e_1\wedge e_4-e_2\wedge
e_3)$. This shows that $\Pi$ is either an $\alpha$ or a
$\beta$-plane depending on whether the basis $(e_1,e_2,e_3,e_4)$
determines the given orientation of $T_pM$ or the opposite one.

\smallskip

 \begin{lemma}\label{ACS J} \rm{\cite{DGMY}}
Let $M$ be a 4-manifold with a neutral metric $g$ and let $X$ and
$Y$ be orthogonal null vector fields which are linearly independent
at every point of $M$. Then the triple $(g,X,Y)$ determines an
orientation and a unique  orientation and $g$-compatible almost
complex structure $J$ on $M$ such that $JX=Y$.
\end{lemma}

In the proof of this lemma, it has been shown that in a
neighbourhood of every point of $M$, there exist vector fields $Z,T$
such that:
\begin{enumerate}
\item[$(i)$] $(X,Y,Z,T)$ is a local frame of the tangent bundle $TM$;

\item[$(ii)$] $g(X,Z)=1$, \quad $g(X,T)=0$,\quad  $g(Y,Z)=0$, \quad $g(Y,T)=1$.
\end{enumerate}
Also, the orientation determined by $(X,Y,Z,T)$ does not depend on
the choice of the vector fields $Z,T$ and we shall say that it is
determined by the triple $(g,X,Y)$. Set $a=g(Z,Z)$, $b=g(T,T)$,
$c=g(Z,T)$ and

\begin{equation}\label{ob}
\begin{array}{lll}
E_1=\displaystyle{\frac{1-a}{2}}X + Z, \quad E_2=\displaystyle{\frac{1-b}{2}}Y +T-cX, \\[8pt]
E_3=-\displaystyle{\frac{1+a}{2}}X + Z, \quad
E_4=-\displaystyle{\frac{1+b}{2}}Y+ T - c X.
\end{array}
\end{equation}
Then $(E_1,E_2,E_3,E_4)$ is an orthogonal frame, positively oriented
with respect to the orientation determined by $(g,X,Y)$ and such
that  $  g(E_1,E_1)=g(E_2,E_2)=1$, $g(E_3,E_3)=g(E_4,E_4)=-1$. The
almost complex structure $J$ for which $JE_1=E_2$, $JE_3=E_4$ has
the required properties \cite{DGMY}.

\smallskip

Now, by (\ref{ob})
$$
X\wedge Y=(E_1-E_3)\wedge (E_2-E_4)=s_1^{+}-s_3^{+}.
$$
Thus, $span\{X,Y\}$ is an $\alpha$-plane. Set
$$
U=X-bY+2T.
$$
This vector field is null and orthogonal to $X$.  Moreover
$$
X\wedge U=-bX\wedge Y+X\wedge (E_2+E_4+bY+2cX)=(E_1-E_3)\wedge
(E_2+E_4)=s_1^{-}+s_3^{-}.
$$
Hence $span\{X,U\}$ is an $\beta$-plane. Next,  $JU$ is  a null
vector field and
$$
\begin{array}{c}
JX\wedge JU=Y\wedge (bX+2JT)=bY\wedge X +Y\wedge
(JE_2+JE_4-bX+2cY)\\[6pt]
=-(E_2-E_4)\wedge (E_1+E_3)=s_1^{-}-s_3^{-}.
\end{array}
$$
Thus $span\{JX,JU\}$ is also a $\beta$-plane. The vectors
$X,JX=Y,U,JU$ are linearly independent at each point since
$$
X\wedge U\wedge JX\wedge JU=-4E_1\wedge E_2\wedge E_3\wedge E_4.
$$
We sum up these observations in the following

\begin{lemma}\label{beta-planes}
In the notation above:

$(i)$ The vectors $\{X,JX=Y,U,JU\}$ are linearly independent  at
every point.

$(ii)$ $g(X,U)=0,\quad g(JX,U)=2,\quad g(X,X)=g(U,U)=0$.

$(iii)$ The null $2$-planes $span\{X,U\}$ and $span\{JX,JU\}$ are
$\beta$-planes.

$(iv)$ Every $\beta$-plane containing $X$ coincides with
$span\{X,U\}$.

$(v)$ Every $\beta$-plane containing $JX$ coincides with
$span\{JX,JU\}$.

$(vi)$ The plane $span\{X,JX\}$ is an $\alpha$-plane and every
$\alpha$-plane containing $X$ coincides with $span\{X,JX\}$.
\end{lemma}

\noindent {\bf Proof}. We have only to prove $(iv)$, $(v)$ and
$(vi)$. Suppose that $\Pi$ is a $2$-plane containing $X$, and let
$\widetilde X=\gamma X+\lambda Y +\mu U +\nu JU$ be a vector such
that $\{X,\widetilde X\}$ is a basis of $\Pi$. We have by (\ref{ob})
$$
\begin{array}{c}
X\wedge JU=X\wedge (Y-bX+2JT)=(1+2c)(X\wedge Y)+X\wedge (JE_2+JE_4)\\[6pt]
=(1+2c)(s_1^{+}-s_3^{+})-(E_1-E_3)\wedge
(E_1+E_3)=(1+2c)(s_1^{+}-s_3^{+})-(s_2^{+}+s_2^{-}).
\end{array}
$$
Then
$$
X\wedge \widetilde X=[\lambda+\nu(1+2c)](s_1^{+}-s_3^{+})-\nu
s_2^{+}+\mu(s_1^{-}+s_3^{-})-\nu s_2^{-}.
$$
It follows that $\Pi$ is a $\beta$-plane if and only if
$\lambda=\nu=0$, i.e. $\widetilde X=\gamma X +\mu U$ where $\mu\neq
0$. Hence $\Pi$ is a $\beta$-plane if and only if $\Pi=span\{X,U\}$.
Also, $\Pi$ is an $\alpha$-plane exactly when $\mu=\nu=0$, i.e.
$\widetilde X=\gamma X+\lambda Y$ where $\lambda\neq 0$ and $Y=JX$.
This proves $(vi)$.

The proof of $(v)$ is similar.

\hfill$\Box$

The following lemma shows that a fixed null-vector in a complex
vector space with a (2,2) signature Hermitian metric determines a
unique para-hyperhermitian structure.

\begin{lemma}\label{S}
Let $V$ be a $4$-dimensional real vector space with a neutral scalar
product $g$, and let  $I$ be a compatible complex structure. Then
for any null vector $X\neq 0$ in $V$ there is a unique
skew-symmetric endomorphism $S$ of $V$ such that $S^2 =Id$,
$IS=-SI$, $SX=X$.
\end{lemma}

\noindent {\bf Proof}. Note first that the complex structure
determined by the triple $(g,X,Y=IX)$ as in  Lemma~~\ref{ACS J}
coincides with $I$.

Now, we apply Lemma \ref{beta-planes}, and define  $S$ to be the
identity on the unique null $\beta$-plane containing $X$, and minus
the identity on the unique null $\beta$-plane containing $IX$. Then
it is easy to check that $S$  has the required properties.

Conversely, suppose that $S$ is an endomorphism with these
properties. Then $S$ is an involution different from $\pm Id$.
Denote the eigenspace of $S$ corresponding to the eigenvalues $+1$
and $-1$ by ${V^{+}}$ and $V^{-}$. Clearly, the complex structure
$I$ is an isomorphism of  ${V^{+}}$ onto $V^{-}$ since it
anti-commutes with $S$. In particular, $dim V^{+}=dim V^{-}=2$. Note
also that the spaces ${V^{+}}$ and $V^{-}$ are isotropic since $S$
is skew-symmetric. Hence $V^{+}$ is either an $\alpha$ or a
$\beta$-plane, and similarly for $V^{-}$. The space $V^{+}$ does not
contain $IX$, hence it is a $\beta$-plane by Lemma~
\ref{beta-planes}. Similarly,  $V^{-}$ is also a $\beta$-plane. This
proves the uniqueness.

\hfill$\Box$

 Using the above algebraic observations, we
provide a  classification of compact neutral Hermitian surfaces
admitting a nowhere vanishing  null vector field.

\begin{theorem}\label{compact-1-field} Let  $(M, g, I)$ be a compact neutral Hermitian surface with
nowhere vanishing null vector field $X$. Then the complex surface
$(M,I)$ is one of the following:

(i) a complex torus;

(ii) a primary Kodaira surface;

(iii) a minimal properly elliptic surface of odd first Betti number;

(iv)  an Inoue surface of type $S^0$ or $S^{\pm}$ without curves;

(v) a Hopf surface;

\smallskip

Conversely, each of the smooth manifolds underlying the complex
surfaces (i) - (iii), an Inoue surface of type $S^{+}$ or a primary
Hopf surface carries a neutral Hermitian metric with a nowhere
vanishing null  vector field.

\smallskip

\end{theorem}

\noindent {\bf Proof}.  Setting $T=IS$, we get an almost
para-hyperhermitian structure $(g,I,S,T)$ on $M$. Let $\Omega_S$ and
$\Omega_T$ be the corresponding fundamental $2$-forms defined by
$\Omega_S(X,Y)=g(SX,Y)$ and similarly for $\Omega_T$.  Then the
$2$-form $\Omega_S+i\Omega_T$ is non-vanishing and of type $(2,0)$
with respect to $I$ which implies the vanishing of the first Chern
class of $(M,I)$. Now the classification of compact complex surfaces
with vanishing first Chern class provides the types $(i) - (v)$ in
the theorem and also the $K3$- surfaces, see, e.g., \cite[Theorem
8]{DGMY}. But we exclude the latter from the list since they do not
admit non-vanishing global vector fields due to the fact that their
Euler characteristic does not vanish .

\smallskip

That the manifolds underlying the surfaces (i) - (iii), an Inoue
surface of type $S^{+}$  or a primary Hopf surface admit
 neutral Hermitian metrics with a null vector field follows
from Theorem~\ref{existence of 2 fields} proved in the next section.

\hfill$\Box$

\noindent {\bf Remark}. If $X$ is Killing and $\mathcal{L}_X I=0$,
 then $\mathcal{L}_X S =0$, where $S$ is the
endomorphism of $TM$ defined by means of $X$ as in Lemma~\ref{S},
hence $X$ preserves the almost para-hypercomplex structure. Indeed,
notice first that the $(+1)$-eigenbundle $S^+$ of $S$ is integrable,
see \cite{Cald}, \cite{DW}. Hence there is a vector field $Y$ in a
neighbourhood of each point of $M$ such that $span\{X,Y\}=S^{+}$ and
$[X,Y]=0$. Then $span\{IX,IY\}$ is the $(-1)$-eigenspace of the
involution $S$. Clearly $[X,SX]=S[X,X]=0$ and $[X,SY]=S[X,Y]=0$, so
$({\mathcal L}_{X}S)(Z)=0$ for $Z\in S^{+}$.  By assumption
$\mathcal{L}_X I = 0$, i.e. $[X,IZ]=I[X,Z]$ for every $Z$, so
$[X,IX]=I[X,X] = 0$. Then $[X,SIX]=-[X,IX]=0$ and
$[S,SIY]=-[X,IY]=-I[X,Y]=0$, hence $({\mathcal L}_{X}S)(Z)=0$ for
$Z\in S^{-}$.

\medskip

It is well known that on a compact K\"ahler manifold every Killing
vector field is real holomorphic.  On general Hermitian manifolds
existence of a Killing vector field which is not real holomorphic,
has a significant consequences for the curvature. We adapt here the
relevant considerations for the split-signature 4-dimensional case.

\begin{theorem}\label{RH-ASD}
Let $(M, g, I)$ be a neutral Hermitian surface and $X$ a complete
null Killing vector field, which is not  real holomorphic. Then the
metric $g$ is anti-self-dual.
\end{theorem}

\noindent {\bf Proof}. Let ${\varphi}_t$ be the global flow of
isometries of $M$ generated by $X$.  Then
$I_t=(\varphi_t)_{\ast}\circ I\circ (\varphi_{-t})_{\ast}$, $t\in
{\mathbb R}$, are compatible complex structures on $M$. Suppose that
the Lie derivative ${\mathcal L}_{X}I\neq 0$ at a point $p\in M$.
Then ${\mathcal L}_{X}I\neq 0$ in a neighborhood $U$ of $p$. For a
point $x\in U$ we can take three different real numbers
$t_1,t_2,t_3$ closed enough to {\color{blue}$0$} such that
$I_{t_i}(x)\neq \pm I_{t_j}(x)$. It follows that the self-dual
component $W_{+}$ of the Weyl tensor vanishes at the point $x$. This
fact has been proved in \cite{S91} in the positive-definite case,
but the same proof using spinors works in the split signature case
too. Hence $W_{+}$ vanishes on $U$ and by the proof of
\cite[Corollary 1.6]{Pon}, $W_{+}$ vanishes identically. Note that
this result is based on \cite[Proposition 1.3]{Pon} which, although
stated in the positive definite case, holds true in the split
signature case with the same proof  replacing the usual twistor
space by the so-called hyperbolic twistor space \cite{BDM}).

\hfill$\Box$

\noindent {\bf Remark}. As is well-known \cite{O'Neill}, on a
complete pseudo-Riemannian manifold, every Killing vector field is
complete.

\section{ Neutral $4$-manifolds with two null conformal Killing vector fields}

\medskip

Recall that a vector field  $K$ on  a pseudo-Riemannian manifold
$(M,g)$  is called conformal Killing if if it satisfies
$$g(\nabla_A K, B) + g(\nabla_B K, A) = \frac{2}{dim\,M} div(K)g(A,B),
\quad A,B\in TM,$$  where $\nabla$ is the Levi-Civita connection.

The following lemma is an extension of a result of Dunajsky-West
\cite[Lemma 1]{DW} and  D. Calderbank \cite[Lemma 4.1]{Cald}.

\begin{lemma}\label{David} Let $K$ and $L$ be two orthogonal null
vector fields that are linearly independent at every point.

$(i)$ If $K$ is conformal Killing, then $\nabla_{A} B \in span\{K,
L\}$ for $A,B \in span\{K, L\}$.

\smallskip

\smallskip

Suppose that both $K$ and $L$ are conformal Killing vector fields.
Then:

\smallskip

 $(ii)$ If $K$ and $L$ commute, the distribution $span\{K,L\}$ is
 parallel and $\nabla_{L}K=\nabla_{K}L=0$.

 $(iii)$ If the distribution $span\{K,L\}$ is
 parallel, the vector fields $K$ and $L$ are Killing and commute.
\end{lemma}

\noindent {\bf Proof}. $(i)$ At every point, the null-plane
$\Pi=span\{K, L\}$ is a Lagrangian subspace, so $\Pi=\Pi^{\perp}$.
Hence $A\in\Pi$ if and only if $A\perp\Pi$. We have
$$
g(\nabla_{K}K,K)=g(\nabla_{K}L,L)=g(\nabla_{L}K,K)=g(\nabla_{L}L,L)=0
$$
because $g(K,K)=g(L,L)=0$. Moreover,
$$
0=g(K,\nabla_{L}K)=-g(\nabla_{K}K,L)=g(K,\nabla_{K}L)
$$
since $K$ is conformal Killing and orthogonal to $ L$. Thus,
$\nabla_{K}K$ and $\nabla_{K}L$ are orthogonal to $span\{K, L\}$.
Moreover,  we have
$$
g(\nabla_{L}K,L)=-g(\nabla_{L}K,L),
$$
hence
$$
g(\nabla_{L}L,K)=-g(L,\nabla_{L}K)=0.
$$
Therefore $\nabla_{L}K$ and $\nabla_{L}L$ are orthogonal to
$span\{K, L\}$.

\medskip

\noindent $(ii)$ For every tangent vector $Z$
$$
g(\nabla_{Z}K,L)=-g(K,\nabla_{Z}L)=g(\nabla_{K}L,Z)=g(\nabla_{L}K,Z)=-g(\nabla_{Z}K,L)
$$
Hence
$$
g(\nabla_{Z}K,L)=0.
$$
Therefore $\nabla_{Z}K\in span\{K,L\}$. Similarly, $\nabla_{Z}L\in
span\{K,L\}$. Moreover,
$$
g(\nabla_{K}L,Z)=-g(\nabla_{Z}L,K)=g(L,\nabla_{Z}K)=0.
$$
Hence $\nabla_{K}L=0$ and $\nabla_{L}K=0$.

\medskip

\noindent $(iii)$ If the distribution $span\{K,L\}$ is parallel,
then for every $Z\in TM$
$$
g(\nabla_{L}K,Z)=-g(\nabla_{Z}K,L)=0.
$$
Thus $\nabla_{L}K=0$. Similarly $\nabla_{K}L=0$. Hence $[K,L]=0$.
Also,
$$
\begin{array}{l}
0=g(\nabla_{Z}K,L)=-g(\nabla_{Z}L,K)=-\frac{1}{2}
div(L)g(Z,K)+g(\nabla_{K}L,Z)\\[6pt]
 =-\frac{1}{2}
div(L)g(Z,K)+ \frac{1}{2} div(K)g(Z,L)-g(\nabla_{Z}K,L)\\[6pt]
= -\frac{1}{2} div(L)g(Z,K)+ \frac{1}{2} div(K)g(Z,L).
\end{array}
$$
Since  $K,L$ are linearly independent at each point, it follows that
$div(L)=div(K)=0$. Thus, $K$ and $L$ are Killing.

\hfill$\Box$

\smallskip

It follows from Lemmas~\ref{ACS J} - \ref{S}  that a neutral
4-manifold with two orthogonal, pointwise linearly independent, null
vector fields admits an almost para-hyperhermitian structure. When
the vector fields are conformal Killing we have a stronger result.

\begin{theorem}\label{thm1}
Let $(M,g)$ be a neutral  4-manifold with two orthogonal, pointwise
linearly independent, null conformal Killing vector fields $X$ and
$Y$. Then $M$ admits a para-hyperhermitian structure $(g,I,S,T)$,
where  $I$ is the structure determined by $\{g,X,Y\}$ as in
Lemma~\ref{ACS J}.

Moreover, if  the vector fields $X$ and $Y$ commute, then they are
real holomorphic and Killing.
\end{theorem}

\noindent {\bf Proof}.  Consider $M$ with the orientation determined
by $(g,X,Y)$ and let $I$ be the compatible almost complex structure
from Lemma~\ref{ACS J}, so $Y = IX$. By Lemma~\ref{beta-planes}, in
the vicinity of every point, there is a null vector field $U$
orthogonal to $X$ such that the vector fields $X,IX,U,IU$ have the
properties $(i)$ - $(v)$ of this lemma. Then the $(1,0)$-vector
fields $X - iIX$ and $U-iIU$ constitute a frame for $(1,0)$-vectors.
Hence, for the integrability of $I$,  it is enough to check that
$g(\nabla_Z V, W) = 0$ for any choice of $ Z,V,W$ among $2X^{1,0}=X
- iIX$ and $2U^{1,0}=U-iIU$. Now, since (1,0)-vectors are isotropic,
for any (1,0)-vector fields $Z,U,W$, we have the following:

\vspace{.1in}

1. $g(\nabla_Z V, W) = -g(V, \nabla_Z  W)$, and, as a consequence,
$g(\nabla_Z V, V) = 0$.

\vspace{.1in}

\noindent Since $X, IX$ are conformal Killing, it easy to see that

\vspace{.1in}

2. $g(\nabla_Z (X- iIX), V) = - g(\nabla_V (X-iIX), Z)$, and, as a
consequence, $g(\nabla_Z(X-iIX), Z) = 0$.

\vspace{.1in}

From here we get the property $g(\nabla_Z V, W) = 0$ for any choice
of $Z,V,W$ among $ 2X^{1,0} = X - iIX$ and $2U^{1,0} = U-iIU$ by a
case by case argument. The details are as follows.

\smallskip

a) $g(\nabla_{U^{1,0}} X^{1,0}, X^{1,0}) = 0$ (by 1.) and
$g(\nabla_{U^{1,0}} X^{1,0}, U^{1,0}) = 0$ (by 2.)

\vspace{.1in}

b) $g(\nabla_{X^{1,0}} U^{1,0}, X^{1,0}) = 0$ (by 2.) and
$g(\nabla_{X^{1,0}} U^{1,0}, U^{1,0}) = 0$ (by 1.)

\vspace{.1in}

c) $g(\nabla_{X^{1,0}} X^{1,0}, X^{1,0}) = 0$ (by 1.) and
$g(\nabla_{X^{1,0}} X^{1,0}, U^{1,0}) = 0$ (by 1. and then b) )

\vspace{.1in}

d) $g(\nabla_{U^{1,0}} U^{1,0}, U^{1,0}) = 0$ (by 1.) and
$g(\nabla_{U^{1,0}} U^{1,0}, X^{1,0}) = 0$ (by 1. and a) )

\vspace{.1in}

Now, we  define an endomorphism $S$ of $TM$ as $S = Id$ on
$span\{X,U\}$ and $S = -Id$ on $span\{IX,IU\}$. Then
Lemma~\ref{David} implies that the distributions $span\{X,U\}$ and
$span\{IX=Y,IU\}$ are integrable. It follows that $S$ is integrable.
At the end, we set $T=IS$. Then $T$ is integrable since $I$ and $S$
are integrable.

\smallskip

Suppose that the vector fields $X$ and $Y=IX$ commute. Then, by
Lemma~\ref{David}, the vector fields $X$ and $Y$ are Killing. If
$\nabla$ is the Levi-Civita connection of $(M,g)$, we have
$\nabla_XIX=\nabla_{IX}X$. It follows that
$$
g(\nabla_{X}Z,V)=g(\nabla_{IX}Z,V)=0
$$
for $Z,V\in span\{X,IX\}$. Moreover

 $$g(\nabla_XIX,U) =
-g(\nabla_UIX,X)=g(\nabla_UX,IX)$$
 and
 $$
 g(\nabla_XIX,U)=g(\nabla_{IX}X,U)=-g(\nabla_UX, IX)+div(X)
 $$
Hence
$$
g(\nabla_{X}IX,U)=\frac{1}{2}div(X)
$$
Then
$$
\begin{array}{c}
g([U,X],IX)=g(\nabla_{U}X,IX)-g(\nabla_{X}U,IX)=-g(\nabla_{U}IX,X)-g(\nabla_{X}U,IX)\\[6pt]
=2g(\nabla_{X}IX,U)=div(X)
\end{array}
$$
since the vector field $IX$ is conformal Killing, $g(U,X)=0$ and
$g(U,IX)=const$. Moreover
$$
\begin{array}{c}
g(I[IU,X),IX)=g(\nabla_{IU}X,X)-g(\nabla_{X}IU),X)=g(\nabla_{X},IU)\\[6pt]
=-g(\nabla_{IU}X,X)+\frac{1}{2}div(X)g(IU,X)=-div(X)
\end{array}
$$
Therefore
$$
g([U,X]+I[IU,X],IX)=0.
$$
Next,
$$
g([U,X],X)=g(\nabla_{U}X,X)+g(U,\nabla_{X}X)=g(\nabla_{U}X,X)-g(\nabla_{U}X,X)=0.
$$
Also,
$$
\begin{array}{c}
g(I[IU,X])=-g(\nabla_{IU}X,IX)+g(\nabla_{X}IU,IX)=g(\nabla_{IX}X,IU)-g(\nabla_{X}IX,IU)\\[6pt]
=g([X,IX],IU)=0.
\end{array}
$$
Hence
$$
g([U,X]+I[IU,X],X)=0.
$$
It is easy to check that
$$
g([U,X]+I[IU,X],U)=g([U,X]+I[IU,X],IU)=0.
$$
It follows that ${\mathcal L}_{X}I={\mathcal L}_{IX}I=0$.
\hfill$\Box$


\smallskip

Now,  we determine the compact 4-manifolds admitting neutral metrics
with two independent null conformal Killing vector fields.

\smallskip

\begin{theorem}\label{existence of 2 fields}
{\rm{(A)}.}  Let  $(M,g)$ be a compact neutral 4-manifold with two
orthogonal null conformal  Killing vector fields $X$ and $Y$ which
are linearly independent at each point. Then
 $M$ is the underlying smooth
manifold of one of the following complex surfaces:

(i) a complex torus

(ii) a primary Kodaira surface

(iii) a minimal properly elliptic surface of odd first Betti number.

(iv)  an Inoue surface of type $S^0$ or $S^{\pm}$ without curves

(v) a Hopf surface

If in addition the  vector fields $X$ and $Y$ commute, then in
$(iv)$ only $S^+$ is allowed.

\smallskip

\noindent {\rm{(B)}.} Conversely, the smooth manifolds underlying
the complex surfaces $(i)$ - $(iii)$, an Inoue surface of type $S^+$
or a primary Hopf surface admit neutral  metrics with two orthogonal
null Killing vector fields which are linearly independent at each
point.

\end{theorem}

\noindent{\bf Proof}.  (A). By Theorem~\ref{thm1}, $M$ carries a
para-hyperhermitian structure $(g,I,S,T)$ such that $Y=IX$. It
follows from Theorem~\ref{compact-1-field}, applied to the compact
Hermitian surface $(M,I,g)$, that $M$ is one of the manifolds $(i)$
- $(v)$.

The claim for the commuting Killing fields follows from
Theorem~\ref{thm1} and the fact that in $(iv)$ only the Inoue
surfaces of type $S^+$ admit holomorphic vector fields \cite{Inoue}.

\smallskip

(B). The proof for the cases $(i)$ - $(iii)$ follows from the proof
of the second part of \cite[Theorem 7]{DGMY}. Here are the details.

$(i)$.  Let $M$ be the quotient of $\mathbb{C}^2$ by a lattice
$<a_1,a_2,a_3,a_4>$. Take a smooth function $\alpha$ on ${\Bbb C}$
such that $\alpha(z+a_3)=\alpha(z)$, $\alpha(z+a_4)=\alpha(z)$. Then
the metric $$\widetilde g=\alpha\, dzd\overline{z} +
2Re(dzd\overline{w})$$ on $\mathbb{C}^2$ descends to a a split
signature  Ricci flat, K\"ahler  metric $g$ on the torus $M$
\cite{P}. For this metric, the global vector field $W$ on $M$ given
in the standard local coordinates $(z,w)$ on the torus by
$\displaystyle{\frac{\partial}{\partial w}}$ is holomorphic,
parallel and null. Hence  if $K=Re W$, the vector fields $K$ and
$IK=Im W$ have the properties stated in lemma.

\smallskip
$(ii)$.  A primary Kodaira surface \cite{Kodaira}  can be obtained
in the following way. Consider the affine transformations
$\varphi_k(z,w)$ of ${\Bbb C}^2$ given by
$$\varphi_k(z,w) = (z+a_k,w+\overline{a}_kz+b_k),$$ where $a_k$, $b _k$, $k=1,2,3,4$, are complex
numbers such that $$a_1=a_2=0, \quad Im(a_3{\overline a}_4) =
b_1\neq 0,\quad b_2\neq 0.$$ They generate a group $G$ of affine
transformations acting freely and properly discontinuously on ${\Bbb
C}^2$ and $M$ is the  quotient space ${\Bbb C}^2/G$ for a suitable
choice of $a_k$ and $b _k$ \cite[p.786]{Kodaira}. Taking into
account the identities
$$\varphi_{k\ast}(\frac{\partial}{\partial
z})=\frac{\partial}{\partial z} + \overline{a_k}
\frac{\partial}{\partial w},\quad
\varphi_{k\ast}(\frac{\partial}{\partial
w})=\frac{\partial}{\partial w}$$ we see that every holomorphic
vector field on $M$ is proportional to the vector field $W$ given in
the local coordinates $(z,w)$ as
$\displaystyle{\frac{\partial}{\partial w}}$. As in \cite{P}, set
$\alpha(z)=f(z) -\gamma\,z-\overline{\gamma}\,\overline{z}$, where
$f(z)$ is a smooth function on ${\Bbb C}$ satisfying the identities
$f(z+a_3)=f(z)$, $f(z+a_4)=f(z)$. Then the metric $\widetilde
g=\alpha\, dzd\overline{z} + 2Re(dzd\overline{w})$ on $\mathbb{C}^2$
descends to a split signature K\"ahler, Ricci flat metric $g$ on $M$
for which the holomorphic vector field $W$ is parallel and null (the
metric $g$ is flat if $f\equiv const$). Then the vector fields $K=Re
W$ and $IK=Im W$ have the required properties.

\smallskip

$(iii)$. Suppose $M$ is a minimal properly elliptic surface of odd
first Betti number. Let $\widetilde{SL(2,R)}$ be the universal
covering group of $SL(2,R)$. By \cite[Theorem 7.4]{Wall}, $M$ is the
quotient of $\widetilde{SL(2,R)}\times {\mathbb R}$ by a co-compact
lattice. Let $A',B',C'$ be the frame of left-invariant vector fields
on  $SL(2,R)$  determined by  the matrices
$A'=1/2\left(\begin{array}{ll}
 0& 1\\
 -1&
0
\end{array}\right)$,
$B'=1/2\left(\begin{array}{ll}
 0& 1\\
 1& 0
\end{array}\right)$,
$C'=1/2\left(\begin{array}{ll}
 1& 0\\
 0& -1
\end{array}\right)$.
Denote the lifts of $A',B',C'$ to vector fields on
$\widetilde{SL(2,R)}$ by $A,B,C$. These vector fields satisfy the
following commutation relations
$$
[A,B]=C, \quad [B,C]=-A,\quad [C,A]=B.
$$
Hence, if $(\alpha,\beta,\gamma)$ is the dual frame of $(A,B,C)$,
$$
d\alpha=\beta\wedge\gamma,\quad d\beta=\alpha\wedge\gamma,\quad
d\gamma=-\alpha\wedge\beta.
$$
Set $V=\displaystyle{\frac{d}{dx}}$, the standard vector field on
${\mathbb R}$, and  let $I$ be the invariant almost complex
structure on $\widetilde{SL(2,R)}\times {\mathbb R}$ given by $IV=C,
IA=B$. If $\theta$ is dual to $V$, $d\theta=0$ and  the basis of
$(1,0)$-forms for the invariant almost complex structure $I$ is
given by $\theta+i\alpha, \beta+i\gamma$.  Notice also that
$$
d(\beta+i\gamma)= -i\alpha\wedge(\beta+i\gamma), \quad
d(\theta+i\alpha) = \beta\wedge (\beta+i\gamma),
$$
so $I$ is integrable. Now consider the form
$$\Omega=(\theta+i\alpha)\wedge(\beta+i\gamma)$$ This is a
$(2,0)$-form and from above it follows
$d\Omega=-\theta\wedge\Omega$. Note that $\omega =
\theta\wedge\alpha-\beta\wedge\gamma$  satisfies $d\omega =
-\theta\wedge\omega$, so the pair $(\Omega,\omega)$ defines a
para-hyperhermitian structure by Lemma \ref{APHHS}, which is left
invariant (and locally conformally  para-hyperk\"ahler). In fact the
metric defined by $\omega$ and $I$ is the bi-invariant metric $g =
\theta^2+\alpha^2-\beta^2-\gamma^2$. One can  directly check that
$$\mathcal{L}_Vg = \mathcal{L}_Ag = \mathcal{L}_Bg = \mathcal{L}_Cg
= 0$$ which also follows from the fact that $V,A,B,C$ are left
invariant. Then the neutral Hermitian structure $(g,I)$ and the
vector fields $K=V+B$ and $IK=A+C$ descend to $M$ giving a neutral
Hermitian structure with two orthogonal null Killing vector fields
which are linearly independent at each point.

\smallskip

$(iv)$.  In order to discuss the Inoue surfaces $S^{+}$,  we recall
first their definition \cite{Inoue}. Take a matrix $N=(n_{ij})\in
GL(2,{\Bbb Z})$ with $det\,N=1$ having two real eigenvalues
$\alpha>1$ and $\alpha^{-1}$. Note that $\alpha$ is a irrational
number. Choose real eigenvectors $(a_1,a_2)$ and $(b_1,b_2)$
corresponding to $\alpha$ and $\alpha^{-1}$, respectively. Take
integers $p,q,r$, $r\neq 0$ and a complex number $t$. Let
$(c_1,c_2)$ be the solution of the equation
\begin{equation}\label{cc}
\varepsilon
(c_1,c_2)=(c_1,c_2)N^{tr}+(e_1,e_2)+\frac{b_1a_2-b_2a_1}{r}(p,q)
\end{equation}
where $N^{tr}$ is the transpose matrix of $N$ and
$$
e_k=\frac{1}{2}n_{k1}(n_{k1}-1)a_1b_1+\frac{1}{2}n_{k2}(n_{k2}-1)a_2b_2+n_{k1}n_{k2}b_1a_2,
\quad k=1,2.
$$
Let $G^{+}=G^{+}_{N,p,q,r;t}$ be the group generated by the
following automorphisms of ${\Bbb C}\times {\bf H}$, ${\bf H}$ being
the upper half-plane:
\begin{equation}\label{G}
\begin{array}{l}
g_0(z,w)=(\varepsilon z+\frac{1}{2}(1+\varepsilon)t,\alpha w)\\[6pt]
g_k(z,w)=(z+b_kw+c_k,w+a_k),\, k=1,2,\\[6pt]
g_3(z,w)=(z+\displaystyle{\frac{b_1a_2-b_2a_1}{r}},w).
\end{array}
\end{equation}
The group $G^{+}$ acts properly discontinuously and without fixed
points in view of (\ref{cc}) and the fact that $(a_1,b_1)$ and
$(a_2,b_2)$ are linearly independent vectors. The quotient
$S^{+}_{N,p,q,r;t}=({\Bbb C}\times {\bf H})/G^{+}$ is a compact
complex surface, known as an Inoue surface of type $S^{+}$.

Set $t_2=Im\,t$ and
\begin{equation}\label{s+}
\alpha_1=dx-\frac{1}{v}(y-t_2\frac{\ln v}{\ln\alpha})du,\quad
\alpha_2=dy-\frac{1}{v}(y-t_2\frac{\ln v}{\ln\alpha})dv,\quad
\alpha_3=\frac{du}{v}, \quad \alpha_4=\frac{dv}{v},
\end{equation}
where $z=x+iy$ and $w=u+iv$. These forms are linearly independent
and invariant under the action of the group $G^{+}$. Note also that
$\alpha_1+i\alpha_2$ and $\alpha_2+i\alpha_4$ constitute a frame of
$(1,0)$-forms. Moreover,
$$
d\alpha_1=\alpha_3\wedge\alpha_2-\frac{t_2}{\ln\alpha}\alpha_3\wedge
\alpha_4,\quad d\alpha_2=\alpha_4\wedge\alpha_2,\quad
d\alpha_3=\alpha_3\wedge\alpha_4,\quad d\alpha_4=0.
$$
Set
$$
\Omega_1=\alpha_1\wedge\alpha_3+\alpha_2\wedge\alpha_4,\quad
\Omega_2=\alpha_1\wedge\alpha_3-\alpha_2\wedge\alpha_4,\quad
\Omega_3=\alpha_1\wedge\alpha_4+\alpha_2\wedge\alpha_3.
$$
Then
$$
-\Omega_1^2 = \Omega_2^2 =
\Omega_3^2=2\alpha_1\wedge\alpha_2\wedge\alpha_3\wedge\alpha_4,
\quad \Omega_l \wedge\Omega_m = 0, \,1\leq l \neq m\leq 3, \quad
d\Omega_l=-\alpha_4\wedge\Omega_l.
$$
Therefore, by Proposition~\ref{phe},  $\Omega_1,\Omega_2,\Omega_3$
define an $G^{+}$-invariant para-hyperhermitian structure on ${\Bbb
C}\times {\bf H}$ with the neutral metric $g(X,Y)=\Omega_1(X,JY)$
which is locally conformally para-hyperk\"ahler since its Lie form
$\theta=-\alpha_4$ is closed. This structure descends to a
para-hyperhermitian structure on the Inoue surface $S^{+}$.

\smallskip

The frame $(X_1,...,X_4)$ of $G^{+}$-invariant vector fields dual to
the frame $(\alpha_1,...,\alpha_4)$ is given by
$$
\displaystyle{X_1=\frac{\partial}{\partial x},\>
X_2=\frac{\partial}{\partial y},\> X_3=(y-t_2\frac{\ln
v}{\ln\alpha})\frac{\partial}{\partial x}+v\frac{\partial }{\partial
u},\> X_4=(y-t_2\frac{\ln v}{\ln\alpha})\frac{\partial}{\partial
y}+v\frac{\partial }{\partial v}}
$$
The non-zero Lie brackets of these vector fields are
$$
[X_2,X_3]=X_1,\quad [X_2,X_4] = X_2, \quad  [X_3,X_4] = -X_3.
$$

If $I$ is the standard complex structure on ${\Bbb C}\times {\bf
H}$,  clearly $IX_1=X_2$, $IX_3=X_4$.  Set $K=X_1$. Then
$\mathcal{L}_K I = \mathcal{L}_K \Omega_1 = 0$ , so $K$ is Killing
and real holomorphic vector field. One can check that $IK=X_2$ is
also Killing and real holomorphic with $[K,IK] = 0$.

\smallskip

  For any non-negative smooth function $f$ on $M$ with $X_1(f)=X_2(f)=0$ consider
the non-degenerate real $(1,1)$-form
$\omega_f=\Omega_1+f\alpha_3\wedge \alpha_4$. It defines a neutral
metric compatible with $I$ and one can check as above that the
vector fields $K,IK$ are null and Killing with respect to this
metric. It is easy to check that the (2,0)-form $\Omega =
(\alpha_1+i\alpha_2)\wedge(\alpha_3+i\alpha_4)$ and the real
$(1,1)$-form $\omega_f$ satisfy the identities
$$\Omega\wedge\overline\Omega= -2\omega_f^2 \ , \
\Omega\wedge\omega_f=0 \ , \ d\Omega=-\alpha_4\wedge\Omega \ , \
d\omega_f = -\alpha_4\wedge\omega_f.$$ Hence, by Lemma \ref{APHHS},
$\Omega$ and $\omega_f$ define a para-hyperhermitian structure on
$M$.

\smallskip

\smallskip

$(v).$ Recall that every primary Hopf surface is diffeomorphic to
$S^1\times S^3$. Consider $S^3$ as the Lie group $SU(2)$ and let
$X_2,X_3,X_4$ be left-invariant vector fields defining its Lie
algebra with the commutator relations
$$
[X_2,X_3]=X_4, \hspace{.2in} [X_3,X_4]=X_2, \hspace{.2in}
[X_4,X_2]=X_3
$$
Set $X_1=\displaystyle{\frac{\partial}{\partial x}}$ where $x$ is
the standard coordinate $e^{it}\to t$ on $S^1$. Denote by $g$ the
metric on $S^1\times S^3$ for which the frame $X_1,...,X_4$ is
orthogonal and $||X_1||^2=||X_2||^2=1$, $||X_3||^2=||X_4||^2=-1$.
This metric is Hermitian with respect to the complex structure
defined by $IX_1=X_2$, $IX_3=X_4$.  It is easy to check that
${\mathcal L}_{X_k}g=0$, $k=1,2,3,4$. Then $X=X_1+X_3$ and
$IX=X_2+X_4$ are two orthogonal null Killing vector fields.

\hfill$\Box$

\smallskip

\noindent {\bf Remark}. We  show that the  Hopf surface  $S^1\times
S^3=S^1\times SU(2)$  does not admit any left-invariant metric with
two non-collinear, left-invariant isotropic and orthogonal Killing
vector fields that commute. To do this we will use the notations
above. Suppose that $g$ is a left-invariant metric (of arbitrary
signature) on $S^1\times SU(2)$  and let $ X=x_0X_0+
x_1X_1+x_2X_2+x_3X_3$ and $Y=y_0X_0+ y_1X_1+y_2X_2+x_3X_3$ be
left-invariant vector fields satisfying the above properties. Then
$$0=[X,Y]=
(x_2y_3-x_3y_2)X_1+(x_3y_1-x_1y_3)X_2+(x_1y_2-x_2y_1)X_3$$ and
therefore the vectors $U=x_1X_1+x_2X_2+x_3X_3$ and
$V=y_1X_1+y_2X_2+x_3X_3$ are collinear. Without loss of generality
we may assume that $V=aU$ and since $X=x_0X_0+U$ and $Y=y_0X_0+aU$
are isotropic and orthogonal we get

\begin{equation}\label{H}
\begin{array}{l}
x_0^2\parallel X_0\parallel^2+2x_0g(X_0,U)+ \parallel U\parallel^2=0\\[6pt]
y_0^2\parallel X_0\parallel^2+2ay_0g(X_0,U)+ a^2\parallel U\parallel^2=0\\[6pt]
x_0y_0\parallel X_0\parallel^2+ (ax_0+y_0)g(X_0,U)+ a\parallel
U\parallel^2=0.
\end{array}
\end{equation}
 Having in mind that $y_0\neq ax_0$ ($X$ and $Y$ are
non-collinear) it follows respectively from the first and second,
and from the first and third equation in (\ref {H}) that

\begin{equation}\label{H1}
\begin{array}{l}
(y_0+ax_0)\parallel X_0\parallel^2 + 2ag(X_0,U)=0\\[6pt]
 x_0\parallel X_0\parallel^2+ g(X_0,U)=0.
\end{array}
\end{equation}

Hence $$(y_0-ax_0)\parallel X_0\parallel^2=0 \Rightarrow \parallel
X_0\parallel^2=0\Rightarrow g(X_0,U)= \parallel U\parallel^2=0 .$$

Note that a left-invariant vector field $A$ is Killing with respect
to the left-invariant metric $g$ if and only if
\begin{equation}\label{H2}
\begin{array}{l}
g([A,B],C)+g([A,C],B)=0
\end{array}
\end{equation}
for all left-invariant vector fields $B$ and $C$. In particular, the
vector field  $X_0$ is Killing since it commutes with all
left-invariant vector fields. Hence $U=x_1X_1+x_2X_2+x_3X_3$ is a
Killing vector field too and substituting $A=U, B=X_0$ in (\ref
{H2}) we get $ g([U,C],X_0)=0$ for $C=X_1,X_2,X_3$. Set
$a_i=g(X_0,X_i), i=1,2,3.$. Then using the commutation relations for
$X_1,X_2,X_3$ we obtain
$$x_2a_1-x_1a_2=x_3a_1-x_1a_3=x_3a_2-x_2a_3=0.$$
These identities together with  $0=g(U,X_0)=x_1a_1+x_2a_2+x_3a_3$
imply that either $x_1=x_2=x_3=0$, i.e. $U=0$ or $a_1=a_2=a_3=0$,
i.e. $X_0=0$ since $g(X_0,X_0)=0$. In both cases it follows that the
vector fields $X$ and $Y$ are collinear which is a contradiction.

\vskip 20pt


\noindent Johann Davidov

\noindent Institute of Mathematics and Informatics

\noindent Bulgarian Academy of Sciences

\noindent 1113 Sofia, Bulgaria

\noindent jtd@math.bas.bg

\medskip

\noindent Gueo Grantcharov

\noindent Department of Mathematics

\noindent Florida International University

\noindent Miami, FL 33199

\noindent grantchg@fiu.edu

\medskip

\noindent Oleg Mushkarov

\noindent Institute of Mathematics and Informatics

\noindent Bulgarian Academy of Sciences

\noindent 1113 Sofia, Bulgaria

\noindent

and

\noindent Sauth-West University

\noindent  2700 Blagoevgrad, Bulgaria

\noindent muskarov@math.bas.bg


\begin{thebibliography}{00}

\bibitem{AG} V. Apostolov, M. Gualtieri, {\it Generalized K\"ahler manifolds, commuting complex structures, and split tangent
bundles},  Comm. Math. Phys. {\bf 271} (2007),  561-575.


\bibitem{Besse} A. Besse, {\it Einstein manifolds}, Classics in Mathematics, Springer-Verlag, Berlin, 2008.

\bibitem{BDM} D. Blair, J. Davidov, O. Mushkarov, {\it Hyperbolic
twistor space}, Rocky Mount. J. Math. {\bf 35} (2005), 1437-1465.

\bibitem{Boyer} C. Boyer, {\it A note on hyperhermitian
foor-manifolds}, Proc. Amer. Math. Soc. {\bf 102} (1988), 157-164.

\bibitem{Cald} D. Calderbank, {\it
Selfdual 4-manifolds, projective surfaces, and the Dunajski-West
construction.} SIGMA Symmetry Integrability Geom. Methods Appl. {\bf
10} (2014), Paper 035, 18 pp.

\bibitem{DGMY09} J. Davidov, G. Grantcharov, O. Mushkarov, M. Yotov,
{\it Para-hyperhermitian surfaces}, Bull. Math. Soc. Sci. Math.
Roumanie {\bf 52} (100) (2009), 281-289.

\bibitem{DGMY}  J. Davidov, G. Grantcharov, O. Muskarov, M. Yotov,
{\it Compact complex surfaces with geometric structures related to
split quaternions}, {\textrm Nucl. Phys. B} 865 (2012), 330--352.

\bibitem{DW} M. Dunajski, S. West, {\it Anti-self-dual conformal structures with null Killing vectors from projective structures}, Comm. Math. Phys. {\bf 272} (2007),  85-118.

\bibitem{DW1}  M. Dunajski, M. West, {\it Anti-self-dual conformal structures in neutral signature,} Recent developments in pseudo-Riemannian geometry, ESI Lect. Math.
Phys., Eur. Math. Soc., Z\"urich, 2008, pp. 113-148


\bibitem{GT} P. Gauduchon, K. P. Tod, {\it Hyper-Hermitian metrics with symmetry}, J. Geom. Phys. {\bf 25} (1998), 291-304.


\bibitem{Hit} N. Hitchin, {\it Hypersymplectic quotiens,} Atti Accad. Sci. Torino cl. Sci. Fis.
Mat. Natur. {\bf 124} (1990), 169-180.


\bibitem{Inoue} M. Inoue, {\it On surfaces of class $VII_0$}, Invent. Math. {\bf 24} (1974), 269-310.

\bibitem{Iv-Zam} S. Ivanov, S. Zamkovoy, {\it Para-hermitian and para-quaternionic
manifolds}, Diff. Geom. Appl. {\bf 23} (2005), 205-234.

\bibitem{Kamada02} H. Kamada, Self-dual K\"ahler metrics of neutral signature on
complex surfaces, PhD thesis, Tohoku University (2002).

\bibitem{Klin} B. Klingler. {\it Chern's conjecture for special affine manifolds.} Ann. of Math. (2) {\bf 186} (2017) (1) 69 - 95.


\bibitem{Kodaira} K. Kodaira, {\it On the structure of compact complex analytic surfaces I}. Amer. J. Math. {\bf 86} (1964),
751-798.

\bibitem{OV1} H. Ooguri, C. Vafa, {\it Geometry of N=2 strings.} Nucl. Phys. B {\bf 361} (1991),
469-518.


\bibitem{OV2} H. Ooguri, C. Vafa, {\it Self-duality and N=2 string magic.} Modern Phys. Lett.
A {\bf 5} (1990), 1389-1398.




\bibitem{O'Neill} B. O'Neill, {\it Semi-Riemannian geometry. With applications to relativity}.
Pure and Applied Mathematics, {\bf 103}, Academic Press, Inc. , New
York, 1983.

\bibitem{Pon} M. Pontecorvo, {\it Complex structures on Riemannian four
manifolds}, Math. Ann. {\bf 309} (1997), 159-177.

\bibitem{P} J. Petean, {\it Indefinite K\"ahler-Einstein metrics on compact complex surfaces},
Comm. Math. Phys. {\bf 189} (1997), 227-235.


\bibitem{S91} S. Salamon, {\it Special structures on
four-manifolds}, Riv. Mat. Iniv. Parma (4) {\bf 17*} (1991),
109-123.


\bibitem{Wall} C. T. C. Wall, {\it Geometric structures on compact complex
analytic surfaces}, Topology {\bf 25} (1986) 119--153.


\bibitem{Ward} R. S. Ward, {\it Integrable and solvable systems, and relations among them.} New developments in the theory and application of solitons. Philos. Trans. Roy. Soc. London Ser. A {\bf 315} (1985), no. 1533, 451-457.



\end{thebibliography}
\end{document}